\documentclass{amsart}

\usepackage{amsfonts,amsmath,amssymb}

\theoremstyle{plain}
\newtheorem{theorem}{Theorem}[section]
\newtheorem{lemma}[theorem]{Lemma}
\newtheorem{corollary}[theorem]{Corollary}
\newtheorem{conjecture}[theorem]{Conjecture}

\theoremstyle{definition}

\numberwithin{equation}{section}

\newcommand{\Z}{\mathbb{Z}}

\font\bb=msbm10 at9.98pt
\newcommand{\semidirect}{\hbox{$\;$\bb\char'156$\;$}}

\begin{document}

\title[Elliptic hypergeometric series]{Summation and
transformation formulas for elliptic hypergeometric series}

\author[S.O. Warnaar]{S.O. Warnaar}

\address{Instituut voor Theoretische Fysica, Universiteit van Amsterdam,
Valckenierstraat 65, 1018 XE Amsterdam, The Netherlands}
\email{warnaar@wins.uva.nl}

\subjclass{Primary 33D15, 33D67, 33E05; Secondary 05A30}

\begin{abstract}
Using matrix inversion and determinant evaluation techniques
we prove several summation and transformation formulas
for terminating, balanced, very-well-poised, elliptic
hypergeometric series. 
\end{abstract}

\maketitle

\section{Introduction}
In the preface to their book ``Basic Hypergeometric Series''~\cite{GR90},
Gasper and Rahman refer to the enchanting nature  
of the theory of $q$-series or basic hypergeometric series as the
highly infectious ``$q$-disease''.
Indeed, from the time of Heine, about a century and a half ago,
till the present day, many researchers have been pursuing
the task of finding $q$-analogues of classical results 
in the theory of special functions,
orthogonal polynomials and hypergeometric series.
It thus seems somewhat surprising that what appears to be the next
natural line of research, replacing ``$q$-analogue'' by ``elliptic analogue'',
has so-far found very few practitioners.

The elliptic (or ``modular'') analogues of hypergeometric series
were introduced by Frenkel and Turaev~\cite{FT97} in their study of elliptic
$6j$-symbols. These $6j$-symbols, which correspond to certain
elliptic solution of the Yang--Baxter equation found by Baxter~\cite{Baxter73} 
and Date et al.~\cite{DJKMO88},
can be expressed in terms of elliptic generalizations of terminating,
balanced, very-well-poised ${_{10}\phi_9}$ series. 
Moreover, the tetrahedral symmetry of the elliptic $6j$ symbols implies an
elliptic analogue of the famous Bailey 
transformation for ${_{10}\phi_9}$ series.
So far, the only follow up on the work of Frenkel and Turaev appears to
be the paper \cite{SZ99} by Spiridonov and Zhedanov,
who presented several contiguous relations for the elliptic analogue of   
${_{10}\phi_9}$ series and who studied an elliptic version of 
Wilson's~\cite{Wilson91} family of biorthogonal rational functions.
As an independent development towards elliptic analogues,
we should also mention the work by Ruijsenaars~\cite{Ruijsenaars97}
and Felder and Varchenko~\cite{FV99}
we studied an elliptic variant of the $q$-gamma function.

The aim of this paper is to prove several further results for elliptic
hypergeometric series. After an introduction to basic and elliptic
hypergeometric series in section~\ref{secqE}, we use 
section~\ref{secMI}
to derive an elliptic matrix inverse. This matrix inverse,
which generalizes a well-known result from the theory of basic series, 
is used repeatedly in section~\ref{secST} to derive an extensive
list of summation and transformation formulas for terminating,
balanced, very-well-poised, elliptic hypergeometric series.
The ``$q$ limits'' of most of these identities correspond to known
results by Gasper and Rahman, Gessel and Stanton, and Chu.
In section~\ref{secC} we establish an elliptic, multivariable
extension of Jackson's ${_8\phi_7}$ sum associated with the C$_n$ root system,
generalizing the basic case due to Schlosser.
Our proof involves an elliptic extension of a
general determinant evaluation by Krattenthaler.
We conclude the paper with a conjectured C$_n$ Bailey transformation
for elliptic hypergeometric series.

\section{Basic hypergeometric series and their elliptic analogues}\label{secqE}
Assume $|q|<1$ and define the $q$-shifted factorial for all integers $n$ by
\begin{equation*}
(a;q)_{\infty}=\prod_{k=0}^{\infty}(1-aq^k) 
\qquad\text{and}\qquad (a;q)_n=\frac{(a;q)_{\infty}}{(aq^n;q)_{\infty}}.
\end{equation*}
Specifically,
\begin{equation*}
(a;q)_n=
\begin{cases}
\prod_{k=0}^{n-1}(1-aq^k) & \text{$n>0$} \\
1 & \text{$n=0$} \\
1/\prod_{k=0}^{-n-1}(1-aq^{n+k})=1/(aq^n;q)_{-n} & \text{$n<0$.}
\end{cases}
\end{equation*}
With the usual condensed notation
$$(a_1,\dots,a_m;q)_n=(a_1;q)_n\cdots(a_m;q)_n$$
we can define an ${_{r+1}\phi_r}$ basic hypergeometric series as~\cite{GR90}
\begin{equation*}
{_{r+1}\phi_r}\Bigl[\genfrac{}{}{0pt}{}{a_1,a_2,\dots,a_{r+1}}
{b_1,\dots,b_r};q,z\Bigr]=
\sum_{k=0}^{\infty}
\frac{(a_1,a_2,\dots,a_{r+1};q)_k}
{(q,b_1,\dots,b_r;q)_k}z^k.
\end{equation*}
Here it is assumed that the $b_i$ are such that none of the terms in the
denominator of the right-hand side vanishes.
When one of the $a_i$ is of the form $q^{-n}$ ($n$ a nonnegative integer)
the infinite sum over $k$ can be replaced by a sum from $0$ to $n$.
In this case the series is said to be terminating.
A ${_{r+1}\phi_r}$ series is called balanced if 
$b_1\dots b_r=qa_1\dots a_{r+1}$ and $z=q$.
A ${_{r+1}\phi_r}$ series is said to be very-well-poised if 
$a_1 q=a_2 b_1=\dots =a_{r+1}b_r$ and $a_2=-a_3=qa_1^{1/2}$.
In this paper we exclusively deal with balanced, very-well poised series
(or rather, their elliptic analogues) and departing from the standard
notation of Gasper and Rahman's book we
use the abbreviation 
\begin{align*}
{_{r+1}W_r}(a_1;a_4,\dots,a_{r+1};q)&=
{_{r+1}\phi_r}\Bigl[\genfrac{}{}{0pt}{}{a_1,qa_1^{1/2},-qa_1^{1/2},a_4,
\dots,a_{r+1}}
{a_1^{1/2},-a_1^{1/2},qa_1/a_4,\dots,qa_1/a_{r+1}};q,q\Bigr] \\
&=\sum_{k=0}^{\infty}
\frac{1-a_1q^{2k}}{1-a_1} \frac{(a_1,a_4,\dots,a_{r+1};q)_k q^k}
{(q,a_1 q/a_4,\dots,a_1 q/a_{r+1};q)_k},
\end{align*}
where we always assume the parameters in the argument of $_{r+1}W_r$ to obey 
the relation $(a_4\dots a_{r+1})^2=a_1^{r-3}q^{r-5}$.

One of the deepest results in the theory of basic hypergeometric series
is Bailey's transformation~\cite{Bailey29}, \cite[Eq. (III.28)]{GR90}
\begin{multline}\label{q109}
{_{10}W_9}(a;b,c,d,e,f,g,q^{-n};q) \\
=\frac{(aq,aq/ef,\lambda q/e,\lambda q/f;q)_n}
{(aq/e,aq/f,\lambda q/ef,\lambda q;q)_n}\:
{_{10}W_9}(\lambda;\lambda b/a,\lambda c/a,\lambda d/a,e,f,g,q^{-n};q),
\end{multline}
where $$bcdefg=a^3q^{n+2}\qquad\text{and}\qquad\lambda=a^2q/bcd.$$
This identity contains many well-known transformation and 
summation theorems for basic series as special cases.
For example, setting $cd=aq$ (so that $\lambda b/a=1$) and then replacing
$e,f,g$ by $c,d,e$ gives Jackson's $q$-analogue of Dougall's $_7F_6$ 
sum~\cite{Jackson21}, \cite[Eq. (II.22)]{GR90}
\begin{equation}\label{q87}
{_8W_7}(a;b,c,d,e,q^{-n};q)=
\frac{(aq,aq/bc,aq/bd,aq/cd;q)_n}{(aq/b,aq/c,aq/d,aq/bcd;q)_n},
\end{equation}
where $$bcde=a^2 q^{n+1}.$$

To introduce the elliptic analogues of basic hypergeometric series
we need the elliptic function
\begin{equation}\label{defE}
E(x)=E(x;p)=(x;p)_{\infty}(p/x;p)_{\infty},
\end{equation}
for $|p|<1$.
Some elementary properties of $E$ are
\begin{equation}\label{rev}
E(x)=-xE(1/x)=E(p/x)
\end{equation}
and the (quasi)periodicity
\begin{equation}\label{periodicity1}
E(x)=(-x)^k p^{\binom{k}{2}}E(xp^k),
\end{equation}
which follows by iterating \eqref{rev}.

Using definition~\eqref{defE} one can define an elliptic analogue 
of the $q$-shifted factorial by
\begin{equation}\label{EP}
(a;q,p)_n=
\begin{cases}
\prod_{k=0}^{n-1}E(aq^k) & \text{$n>0$} \\
1 & \text{$n=0$} \\
1/\prod_{k=0}^{-n-1}E(aq^{n+k})=1/(aq^n;q,p)_{-n} & \text{$n<0$.}
\end{cases}
\end{equation}
Note that $E(x;0)=1-x$ and hence $(a;q,0)_n=(a;q)_n$.
Again we use condensed notation, setting
$$(a_1,\dots,a_m;q,p)_n=(a_1;q,p)_n\cdots(a_m;q,p)_n.$$

Many of the relations satisfied by the $q$-shifted factorials
(see (I.7)--(I.30) of \cite{GR90}) trivially generalize to the elliptic case.
Here we only list those identities needed later.
The proofs merely require manipulation of the definition of $(a;q,p)_n$;
\begin{subequations}\label{el}
\begin{align}
(aq^{-n};q,p)_n&=(q/a;q,p)_n(-a/q)^n q^{-\binom{n}{2}} \\
(aq^{-n};q,p)_k&=(q/a;q,p)_n(a;q,p)_kq^{-nk}/(q^{1-k}/a;q,p)_n \\
(aq^n;q,p)_k&=(aq^k;q,p)_n(a;q,p)_k/(a;q,p)_n=(aq;q,p)_{n+k}/(a;q,p)_n \\
(a;q,p)_{n-k}&=(a;q,p)_n(-q^{1-n}/a)^k q^{\binom{k}{2}}/(q^{1-n}/a;q,p)_k \\
(a;q,p)_{kn}&=(a,aq,\dots,aq^{k-1};q^k,p)_n.
\end{align}
\end{subequations}
Finally we will need the identity
\begin{equation}\label{periodicity3}
(a;q,p)_n=(-a)^{nk}p^{n\binom{k}{2}}q^{k\binom{n}{2}}(ap^k;q,p)_n,
\end{equation}
which follows from \eqref{periodicity1} and \eqref{EP}.

After these preliminaries we come to Frenkel and Turaev's definition 
of balanced, very-well-poised, elliptic (or modular) hypergeometric 
series~\cite{FT97},
\begin{equation}\label{omega}
{_{r+1}\omega_r}(a_1;a_4,\dots,a_{r+1};q,p)=
\sum_{k=0}^{\infty} \frac{E(a_1q^{2k})}{E(a_1)}
\frac{(a_1,a_4,\dots,a_{r+1};q;p)_k q^k}
{(q,a_1 q/a_4,\dots,a_1 q/a_{r+1};q,p)_k},
\end{equation}
where $(a_4\dots a_{r+1})^2=a_1^{r-3}q^{r-5}$.
Following \cite{FT97} we will stay clear of any convergence problems
by demanding terminating series, i.e., one of the $a_i$ $(i=4,\dots,r+1)$
is of the form $q^{-n}$ with $n$ a nonnegative integer.
Remark that by $E(x;p)E(-x;p)=E(x^2;p^2)$
the above ratio of two elliptic $E$-functions can be written as
\begin{equation*}
\frac{(qa_1^{1/2},-qa_1^{1/2};q,p^{1/2})_k}
{(a_1^{1/2},-a_1^{1/2};q,p^{1/2})_k}.
\end{equation*}
Hence in the $p\to 0$ limit we recover the usual definition
of a balanced, very-well-poised, basic hypergeometric series.

An important result of Frenkel and Turaev is the elliptic analogue of
Bailey's transformation \eqref{q109}.
\begin{theorem}
Let $bcdefg=a^3 q^{n+2}$ and $\lambda=a^2q/bcd.$ Then
\begin{multline}\label{E109}
{_{10}\omega_9}(a;b,c,d,e,f,g,q^{-n};q,p) \\
=\frac{(aq,aq/ef,\lambda q/e,\lambda q/f;q,p)_n}
{(aq/e,aq/f,\lambda q/ef,\lambda q;q,p)_n}\:
{_{10}\omega_9}(\lambda;\lambda b/a,\lambda c/a,\lambda d/a,e,f,g,q^{-n};q,p).
\end{multline}
\end{theorem}
Of course we can again specialize $cd=aq$ to arrive at an elliptic 
Jackson sum.
\begin{corollary}
For $a^2 q^{n+1}=bcde$ there holds
\begin{equation}\label{e87}
{_8\omega_7}(a;b,c,d,e,q^{-n};q,p)=
\frac{(aq,aq/bc,aq/bd,aq/cd;q,p)_n}{(aq/b,aq/c,aq/d,aq/bcd;q,p)_n}.
\end{equation}
\end{corollary}

\section{A matrix inverse}\label{secMI}
Before deriving new summation and transformation
formulas for elliptic hypergeometric series we need to
prepare several, mostly elementary, results for the elliptic
function $E$ of equation~\eqref{defE}. This will result in a
matrix inverse that will be at the heart of all results of the
subsequent section.

The obvious starting point is the well-known addition formula~\cite{WW96}
\begin{multline}\label{Esum}
E(ux)E(u/x)E(vy)E(v/y)-E(uy)E(u/y)E(vx)E(v/x)\\
=\frac{v}{x} E(xy)E(x/y)E(uv)E(u/v).
\end{multline}
By iterating this equation one readily derives the following lemma,
which for $p=0$ reduces to a result of Macdonald, first published
by Bhatnagar and Milne~\cite[Thm. 2.27]{BM97}.
\begin{lemma}
For $n$ a nonnegative integer and $a_j,b_j,c_j,d_j$ ($j=0,\dots,n$)
indeterminates there holds
\begin{multline}\label{EMd}
\sum_{k=0}^n  b_k/c_k E(a_kb_k)E(a_k/b_k)E(c_kd_k)E(c_k/d_k) \\
\times 
\prod_{j=0}^{k-1}E(a_jc_j)E(a_j/c_j)E(b_jd_j)E(b_j/d_j)
\prod_{j=k+1}^n E(a_jd_j)E(a_j/d_j)E(b_jc_j)E(b_j/c_j) \\
=\prod_{j=0}^n E(a_jc_j)E(a_j/c_j)E(b_jd_j)E(b_j/d_j)
-\prod_{j=0}^n E(a_jd_j)E(a_j/d_j)E(b_jc_j)E(b_j/c_j).
\end{multline}
\end{lemma}
\begin{proof}
We carry out induction on $n$.
For $n=0$ the lemma is nothing but \eqref{Esum} with $u=a_0$, $v=b_0$,
$x=c_0$ and $y=d_0$. Now write \eqref{EMd} as $L_n=R_n$ and
assume this to hold for $n\leq m-1$. With the abbreviations
\begin{align*}
f_j &= E(a_jb_j)E(a_j/b_j)E(c_jd_j)E(c_j/d_j) \\
g_j &= E(a_jc_j)E(a_j/c_j)E(b_jd_j)E(b_j/d_j) \\
h_j &= E(a_jd_j)E(a_j/d_j)E(b_jc_j)E(b_j/c_j)
\end{align*}
we then have
\begin{align*}
L_m&=h_m L_{m-1}+\frac{b_m}{c_m} f_m \prod_{j=0}^{m-1}g_j \\
&=h_m \prod_{j=0}^{m-1}g_j-\prod_{j=0}^m h_j
+\frac{b_m}{c_m} f_m \prod_{j=0}^{m-1}g_j  \\
&=\prod_{j=0}^m g_j-\prod_{j=0}^m h_j=R_m,
\end{align*}
where in the second line we have used the induction hypothesis
and in the third line the addition formula \eqref{Esum}
in the form $h_m+b_mf_m/c_m=g_m$.
\end{proof}

Making the substitutions
\begin{equation*}
a_j\to (abd^2)^{1/2},~~b_j\to(ab/c^2)^{1/2},~~c_j\to
(ab)^{1/2}r^j,~~d_j\to(a/b)^{1/2}q^j
\end{equation*}
and using the definition of the elliptic analogue of the
$q$-shifted factorial \eqref{EP}
and the relations \eqref{el} it follows that
\begin{multline*}
\sum_{k=0}^n
\frac{E(aq^kr^k)E(bq^{-k}r^k)}{E(a)E(b)}
\frac{(a/c,c/b;q,p)_k(abd,1/d;r,p)_k q^k}
{(cr,abr/c;r,p)_k(q/bd,adq;q,p)_k} \\
= \frac{E(c)E(ab/c)E(ad)E(bd)}
{E(a)E(b)E(cd)E(abd/c)}
\biggl(1-
\frac{(a/c,bq^{-n}/c;q,p)_{n+1}(abd,dr^{-n};r,p)_{n+1}}
{(bdq^{-n},ad;q,p)_{n+1}(r^{-n}/c,ab/c;r,p)_{n+1}}\biggr).
\end{multline*}
For $p=0$ this corresponds to \cite[Eq. (2.7)]{GR90b} of Gasper and Rahman.
Important will be the specialization obtained by choosing $d=r^n$,
\begin{multline}\label{sum1}
\sum_{k=0}^n
\frac{E(aq^kr^k)E(bq^{-k}r^k)}{E(a)E(b)}
\frac{(a/c,c/b,q,p)_k(abr^n,r^{-n};r,p)_k q^k}
{(cr,abr/c;r,p)_k(qr^{-n}/b,aqr^n;q,p)_k} \\
=\frac{E(c)E(ab/c)E(ar^n)E(br^n)}{E(a)E(b)E(cr^n)E(abr^n/c)}.
\end{multline}
We will use this identity in the next section to prove Theorem~\ref{thmr}.
Now it is needed to obtain the following pair of infinite-dimensional,
lower-triangular matrices, that are inverses of each other
\begin{align*}
f_{n,k}&=\frac{(aq^kr^k,q^kr^{-k}/b;q,p)_{n-k}}
{(r,abr^{2k+1};r,p)_{n-k}} \\
f^{-1}_{n,k}&=(-1)^{n-k}q^{\binom{n-k}{2}}
\frac{E(aq^kr^k)E(q^kr^{-k}/b)}{E(aq^nr^n)E(q^nr^{-n}/b)}
\frac{(aq^{k+1}r^n,q^{k+1}r^{-n}/b;q,p)_{n-k}}
{(r,abr^{n+k};r,p)_{n-k}}.
\end{align*}
For $p=0$ this is \cite[Eqs. (4.4) and (4.5)]{Bressoud88},
\cite[Eqs. (3.2) and (3.3)]{Gasper89} and
\cite[Eq. (4.3)]{Krattenthaler96}.
To derive it from \eqref{sum1} we follow \cite[Sec. 3.6]{GR90}
and set $c=1$ followed by the replacements
$n\to n-l$, $k\to k-l$, $a\to aq^lr^l$ and $b\to bq^{-l}r^l$.
By \eqref{el} one then finds the desired orthogonality relation
\begin{equation}\label{ortho}
\sum_{k=l}^n f^{-1}_{n,k}f_{k,l}=\delta_{n,l}
\end{equation}
with $f$ and $f^{-1}$ as given above.
Finally replacing $r\to q^r$ and using \eqref{el} yields the new
inverse pair
\begin{subequations}\label{invp}
\begin{align}
f_{n,k}&=\frac{E(abq^{2rk})}{E(ab)}
\frac{(aq^n;q,p)_{rk}}{(bq^{1-n};q,p)_{rk}}
\frac{(ab,q^{-rn};q^r,p)_k}{(q^r,abq^{rn+r};q^r,p)_k} \: q^{rk} \\
f^{-1}_{n,k}&=\frac{(b;q,p)_{rn}}{(aq;q,p)_{rn}}
\frac{E(aq^{(r+1)k})E(bq^{(r-1)k})}{E(a)E(b)} \\
&\qquad \times
\frac{(a,1/b;q,p)_k}{(q^r,abq^r;q^r,p)_k}
\frac{(abq^{rn},q^{-rn};q^r,p)_k}{(q^{1-rn}/b,aq^{rn+1};q,p)_k} \: q^k.
\notag
\end{align}
\end{subequations}
This last pair of inverse matrices will be used repeatedly in the next section.
We note that it also follows by the (simultaneous) substitutions
$a\to ab$, $b_i\to aq^i$ and $c_i\to q^{ri}$ in the following elliptic 
analogue of a result due to Krattenthaler~\cite[Eq (1.5)]{Krattenthaler96}.
\begin{lemma}\label{LemKr}
Let $a$ and $b_i,c_i$ ($i\in\Z$) be indeterminates (such that
$c_i\neq c_j$ for $i\neq j$ and $ac_ic_j\neq 1$ for $i,j\in\Z)$.
Then \eqref{ortho} holds with 
\begin{equation*}
f_{n,k}=\frac{\prod_{j=k}^{n-1}E(c_kb_j)E(ac_k/b_j)}
{\prod_{j=k+1}^n c_jE(ac_kc_j)E(c_k/c_j)}
\end{equation*}
and
\begin{equation*}
f^{-1}_{n,k}=\frac{E(c_kb_k)E(ac_k/b_k)}
{E(c_nb_n)E(ac_n/b_n)}
\frac{\prod_{j=k+1}^n E(c_nb_j)E(ac_n/b_j)}
{\prod_{j=k}^{n-1}c_jE(ac_nc_j)E(c_n/c_j)}.
\end{equation*}
\end{lemma}
\begin{proof}
Since for $n=l$ \eqref{ortho} clearly holds we may assume $n>l$ in
the following.
Now let $n>0$ in \eqref{EMd} and make the replacements
$n\to n-l$, $k\to k-l$ and $a_j\to a^{1/2} c_l$, $b_j\to a^{1/2}c_{n+l}$,
$c_j\to b_{j+l}/a^{1/2}$, $d_j\to a^{1/2} c_{j+l}$.
Noting that, in particular, $a_0=d_0$ and $b_n=d_n$ so that the
right-hand side of \eqref{EMd} vanishes, and after performing a few 
trivialities, one finds \eqref{ortho} (with $n>l$) with
$f$ and $f^{-1}$ given by Lemma~\ref{LemKr}.
\end{proof}

\section{Summation and transformation formulas}\label{secST}
Our approach to elliptic hypergeometric summation and transformation
formulas is a standard one in the context of basic hypergeometric
series, see e.g.,~\cite{Chu94,Chu95,GS83,GS86,Krattenthaler96,Riordan79}.
Given a pair of infinite-dimensional, lower-triangular matrices
$f$ and $f^{-1}$, i.e., a pair of matrices such that \eqref{ortho} holds,
the following two statements are equivalent
\begin{equation}\label{fab}
\sum_{k=0}^n f_{n,k} a_k =b_n
\end{equation}
and
\begin{equation}\label{dual}
\sum_{k=0}^n f^{-1}_{n,k}b_k=a_n.
\end{equation}

Our first example of how this may be usefully applied
arises by noting that, thanks to \eqref{sum1}, equation \eqref{dual} holds
for the pair of matrices given in \eqref{invp}
with 
\begin{align*}
a_n&=\frac{(bq;q,p)_{rn}(c,ab/c;q^r,p)_n}
{(a;q,p)_{rn}(abq^r/c,cq^r;q^r,p)_n} \\
b_n&=\frac{(a/c,c/b,q,p)_n(q^r,abq^r;q^r,p)_n}
{(cq^r,abq^r/c;q^r,p)_n(a,1/b;q,p)_n}.
\end{align*}
Hence also \eqref{fab} holds leading to
\begin{multline*}
\sum_{k=0}^n \frac{E(abq^{2rk})}{E(ab)}
\frac{(bq,aq^n;q,p)_{rk}}
{(a,bq^{1-n};q,p)_{rk}}
\frac{(ab,c,ab/c,q^{-rn};q^r,p)_k}
{(q^r,abq^r/c,cq^r,abq^{rn+r};q^r,p)_k} \: q^{rk} \\
=\frac{(a/c,c/b,q,p)_n(q^r,abq^r;q^r,p)_n}
{(cq^r,abq^r/c;q^r,p)_n(a,1/b;q,p)_n}.
\end{multline*}
By \eqref{el} we may alternatively write this in hypergeometric notation
as follows.
\begin{theorem}\label{thmr}
For $r$ a positive integer there holds
\begin{multline*}
{_{2r+6}\omega_{2r+5}}
(ab;c,ab/c,bq,bq^2,\dots,bq^r,aq^n,aq^{n+1},\dots,aq^{n+r-1},q^{-rn};q^r,p) \\
=\frac{(a/c,c/b,q,p)_n(q^r,abq^r;q^r,p)_n}
{(cq^r,abq^r/c;q^r,p)_n(a,1/b;q,p)_n}.
\end{multline*}
\end{theorem}
When $r=1$ this corresponds to the specialization $bc=a$
of the elliptic Jackson sum \eqref{e87}.

With Theorem~\ref{thmr} at hand we are prepared for the proof of
the following quadratic transformation.
\begin{theorem}\label{Etrafo}
Let $bcd=aq$ and $ef=a^2q^{2n+1}$.
When $g=a/b$ or $g=a/e$, there holds
\begin{multline*}
\sum_{k=0}^n \frac{E(aq^{3k})}{E(a)}
\frac{(b,c,d;q,p)_k}{(aq^2/b,aq^2/c,aq^2/d;q^2,p)_k}
\frac{(e,f,q^{-2n};q^2,p)_k}{(aq/e,aq/f,aq^{2n+1};q,p)_k}\: q^k \\
=\frac{(aq^2,a^2q^2/bce,a^2q^2/bdeg,agq^2/cd;q^2,p)_n}
{(a^2q^2/beg,agq^2/c,aq^2/d,a^2q^2/bcde;q^2,p)_n} \\
\times
{_{10}\omega_9}(ag/c;a/c,gq^2/c,beg/a,d,f,g,q^{-2n};q^2,p).
\end{multline*}
\end{theorem}
For $p=0$ and $g=a/b$ this identity reduces to a transformation of
Gasper and Rahman~\cite[Eq. (5.14)]{GR90b}.
We should also remark that
the left-hand side does not depend on $g$ so that the two different
cases actually correspond to a transformation of the right-hand side.
Indeed, the equality of the $g=a/b$ and $g=a/e$ instances of the 
right-hand side is an immediate consequence of the $_{10}\omega_9$
transformation \eqref{E109}.
\begin{proof}
Given the previous remark we only need to prove the $g=a/b$ case 
of the theorem.
Starting point is again the pair of inverse matrices \eqref{invp}
in which we set $r=2$.
The crux of the proof is the observation that equation
\eqref{fab} holds, with 
\begin{equation}\label{a}
a_n=\frac{(b/d,bdq;q^2,p)_n}{(adq^2,aq/d;q^2,p)_n} \:
{_{10}\omega_9}(ad;ad/c,c,dq,dq^2,aq/b,abq^{2n},q^{-2n};q^2,p)
\end{equation}
and
\begin{equation}\label{b}
b_n=\frac{(q^2,abq^2,aq/b;q^2,p)_n(a/c,c/d,dq;q,p)_n}
{(a,1/b,bq;q,p)_n(cq^2,adq^2/c,aq/d;q^2,p)_n}.
\end{equation}
Assuming this is true, we immediately recognize \eqref{dual}
as Theorem~\ref{Etrafo} (with $g=a/b$) under the simultaneous
replacements $b\to a/c, c\to c/d, d\to dq, e\to aq/b$ and $f\to abq^{2n}$.

Of course, it remains to show \eqref{fab} with the above pair of
$a_n$ and $b_n$. Writing $a_n=\sum_j a_{n,j}$ in accordance with the
definition of ${_{10}\omega_9}$, we start with the trivialities
\begin{equation}\label{triv}
b_n=\sum_{k=0}^n f_{n,k}a_k=\sum_{k=0}^n\sum_{j=0}^k f_{n,k}a_{k,j}=
\sum_{j=0}^n\sum_{k=j}^n f_{n,k}a_{k,j}=
\sum_{j=0}^n\sum_{k=0}^{n-j} f_{n,k+j}a_{k+j,j}.
\end{equation}
Using the explicit expressions for $f_{n,k}$ and $a_{n,k}$ as well as the
relations in \eqref{el} this becomes
\begin{multline*}
b_n=\sum_{j=0}^n\frac{(dq,aq^n;q,p)_{2j}(abq^2;q^2,p)_{2j}
(ad,c,ad/c,aq/b,q^{-2n};q^2,p)_j(bq^2/d)^j}
{(a,bq^{1-n};q,p)_{2j}(ad;q^2,p)_{2j}
(q^2,adq^2/c,cq^2,aq/d,abq^{2n+2};q^2,p)_j} \\
\times
{_8\omega_7}(abq^{4j};b/d,bdq^{2j+1},aq^{n+2j},aq^{n+2j+1},q^{-2n+2j};q^2,p).
\end{multline*}
By the elliptic Jackson sum \eqref{e87} we can sum the ${_8\omega_7}$, and
after the usual simplifications we find
\begin{multline*}
b_n=\frac{(abq^2,aq/b;q^2,p)_n}{(adq^2,aq/d;q^2,p)_n} 
\frac{(1/d,dq;q,p)_n}{(1/b,bq;q,p)_n} \\
\times {_{10}\omega_9}(ad;c,ad/c,dq,dq^2,aq^n,aq^{n+1},q^{-2n};q^2,p).
\end{multline*}
By the $r=2$ case of Theorem~\ref{thmr} the ${_{10}\omega_9}$ can be summed
yielding the expression for $b_n$ given in \eqref{b}.
\end{proof}

A result very similar to that of Theorem~\ref{Etrafo} is the following 
cubic transformation, which, for $f=a/b$, provides an elliptic analogue of
\cite[Eq.~(3.6)]{GR90b} by Gasper and Rahman.
\begin{theorem}\label{Etrafo2}
Let $bcd=aq$ and $de=a^2q^{3n+1}$. Then for $f=a/b$ or $f=a/e$ there holds
\begin{multline*}
\sum_{k=0}^n \frac{E(aq^{4k})}{E(a)}
\frac{(b,c;q,p)_k}{(aq^3/b,aq^3/c;q^3,p)_k}
\frac{(d;q,p)_{2k}}{(aq/d;q,p)_{2k}}
\frac{(e,q^{-3n};q^3,p)_k}{(aq/e,aq^{3n+1};q,p)_k}\: q^k \\
=\frac{(aq^3,a^2q^3/bce,a^2q^3/bdef,afq^3/cd;q^3,p)_n}
{(a^2q^3/bef,aq^3f/c,aq^3/d,a^2q^3/bcde;q^3,p)_n} \\
\times
{_{10}\omega_9}(af/c;a/c,fq^3/c,bef/a,d,dq,f,q^{-3n};q^3,p).
\end{multline*}
\end{theorem}
Again we note that the two different cases correspond to the ${_{10}\omega_9}$
transformation \eqref{E109} applied to the right-hand side.
\begin{proof}
By the above remark we only need a proof for $f=a/b$.
The claim is now that if we choose $r=3$ in the pair of 
matrices \eqref{invp} then \eqref{fab} holds, with
\begin{equation}\label{a2}
a_n=\frac{(b^2/a;q^3,p)_n}{(a^2q^3/b;q^3,p)_n} \:
{_{10}\omega_9}(a^2/b;ac/b,a/c,aq/b,aq^2/b,aq^3/b,abq^{3n},q^{-3n};q^3,p)
\end{equation}
and
\begin{equation}\label{b2}
b_n=\frac{(q^3,abq^3;q^3,p)_n(b/c,c;q,p)_n(aq/b;q,p)_{2n}}
{(a,1/b;q,p)_n(acq^3/b,aq^3/c;q^3,p)_n(bq;q,p)_{2n}}.
\end{equation}
Clearly, if this is true we are done with the proof since with these
$a_n$ and $b_n$ equation \eqref{dual} corresponds to the $f=b/a$ case
of Theorem~\ref{Etrafo2} with the replacements
$b\to b/c, d\to aq/b$ and $e\to abq^{3n}$.

To establish \eqref{fab} with the above $a_n$ and $b_n$ we follow
the proof of Theorem~\ref{Etrafo}.
That is, we again write $a_n=\sum_j a_{n,j}$ and use \eqref{triv}.
Inserting the expressions for $f_{n,k}$ and $a_{n,k}$ this yields
\begin{multline*}
b_n=\sum_{j=0}^n
\frac{(aq/b,aq^n;q,p)_{3j}(abq^3;q^3,p)_{2j}
(a^2/b,ac/b,a/c,q^{-3n};q^3,p)_j(b^2q^3/a)^j}
{(a,bq^{1-n};q,p)_{3j}(a^2/b;q^3,p)_{2j}
(q^3,aq^3/c,acq^3/b,abq^{3n+3};q^3,p)_j} \\
\times {_8\omega_7}(abq^{6j};b^2/a,aq^{n+3j},aq^{n+3j+1},aq^{n+3j+2},
q^{-3n+3j};q^3,p).
\end{multline*}
The ${_8\omega_7}$ can be summed by \eqref{e87}, and after
some manipulations involving \eqref{el} we arrive at
\begin{multline*}
b_n=\frac{(abq^3;q^3,p)_n}{(a^2q^3/b;q^3,p)_n}
\frac{(b/a;q,p)_n}{(1/b;q,p)_n}
\frac{(aq/b;q,p)_{2n}}{(bq;q,p)_{2n}} \\
\times 
{_{12}\omega_{11}}(a^2/b;ac/b,a/c,aq/b,aq^2/b,aq^3/b,
aq^n,aq^{n+1},aq^{n+2},q^{-3n};q^3,p).
\end{multline*}
According to Theorem~\ref{thmr} with $r=3$ the ${_{12}\omega_{11}}$ can 
be summed to yield \eqref{b2} as claimed.
\end{proof}

Theorems~\ref{Etrafo} and \ref{Etrafo2} imply several other
quadratic and cubic summation and transformation formulas.

The most obvious ones arise when we demand that $g=1$ in Theorem~\ref{Etrafo}
or $f=1$ in Theorem~\ref{Etrafo2}.
\begin{corollary}\label{cor1}
Let $bcd=aq$ and $ef=a^2q^{2n+1}$. When $b=a$ or $e=a$ there holds
\begin{multline*}
\sum_{k=0}^n \frac{E(aq^{3k})}{E(a)}
\frac{(b,c,d;q,p)_k}{(aq^2/b,aq^2/c,aq^2/d;q^2,p)_k}
\frac{(e,f,q^{-2n};q^2,p)_k}{(aq/e,aq/f,aq^{2n+1};q,p)_k}\: q^k \\
=\frac{(aq^2,a^2q^2/bce,a^2q^2/bde,aq^2/cd;q^2,p)_n}
{(a^2q^2/be,aq^2/c,aq^2/d,a^2q^2/bcde;q^2,p)_n}.
\end{multline*}
\end{corollary}
For $p=0$ and $b=a$ this is a summation of Gessel and Stanton
\cite[Eq. (1.4)]{GS83}, and for $p=0$, $e=a$ it corresponds
to \cite[Eq. (1.9), $b\to q^{-2n}$]{Rahman93} by Rahman and
\cite[Eq. (5.1d)]{Chu95} by Chu.
\begin{corollary}
Let $bcd=aq$ and $de=a^2q^{3n+1}$. When $b=a$ or $e=a$ there holds
\begin{multline*}
\sum_{k=0}^n \frac{E(aq^{4k})}{E(a)}
\frac{(b,c;q,p)_k}{(aq^3/b,aq^3/c;q^3,p)_k}
\frac{(d;q,p)_{2k}}{(aq/d;q,p)_{2k}}
\frac{(e,q^{-3n};q^3,p)_k}{(aq/e,aq^{3n+1};q,p)_k}\: q^k \\
=\frac{(aq^3,a^2q^3/bce,a^2q^3/bde,aq^3/cd;q^3,p)_n}
{(a^2q^3/be,aq^3/c,aq^3/d,a^2q^3/bcde;q^3,p)_n}.
\end{multline*}
\end{corollary}
When $b=a$ this is the elliptic analogue of 
\cite[Eq. (5.22), $c\to q^{-3n}$]{Gasper89} of Gasper.
Theorem~\ref{Etrafo2} also leads to a summation formula if we choose
$d=a$. Indeed, 
${_{10}\omega_9}(af/c;a/c,fq^3/c,bef/a,a,aq,f,q^{-3n};q^3,p)$
with $f=a/b$, $bc=q$ and $e=aq^{3n+1}$ becomes
${_7\omega_6}(a^2/q;a/b,ab/q,aq^{3n+1},q^{-3n};q^3,p)$
which, by \eqref{e87}, evaluates to
\begin{equation*}
\frac{(q^3,a^2q^2,bq,q^2/b;q^3,p)_n}
{(aq,q^2/a,aq^3/b,abq^2;q^3,p)_n}.
\end{equation*}
We therefore conclude the following result, which for $p=0$ corresponds
to \cite[Eq.~(3.7)]{GR90b}.
\begin{corollary}
For $bc=q$ and $e=aq^{3n+1}$,
\begin{multline*}
\sum_{k=0}^n \frac{E(aq^{4k})}{E(a)}
\frac{(b,c;q,p)_k}{(aq^3/b,aq^3/c;q^3,p)_k}
\frac{(a;q,p)_{2k}}{(q;q,p)_{2k}}
\frac{(e,q^{-3n};q^3,p)_k}{(aq/e,aq^{3n+1};q,p)_k}\: q^k \\
=\frac{(aq^2,aq^3,bq,cq;q^3,p)_n}
{(q,q^2,aq^3/b,aq^3/c;q^3,p)_n}.
\end{multline*}
\end{corollary}

The next three results, stated as separate theorems, are somewhat less trivial
as their proof deviates from the standard polynomial argument applicable
in the $p=0$ case.
\begin{theorem}\label{Etrafo3}
For $bcd=a^2q$ and $ef=aq^{n+1}$,
\begin{multline*}
\sum_{k=0}^n\frac{E(aq^{3k})}{E(a)}
\frac{(b,c,d;q^2,p)_k}{(aq/b,aq/c,aq/d;q,p)_k}
\frac{(e,f,q^{-n};q,p)_k}{(aq^2/e,aq^2/f,aq^{n+2};q^2,p)_k}\: q^k \\
\frac{(aq,aq/bc;q,p)_n(aq^{1-n}/b,aq^{1-n}/c;q^2,p)_n}
{(aq/b,aq/c;q,p)_n(aq^{1-n},aq^{1-n}/bc;q^2,p)_n}\\
\times
{_{10}\omega_9}(a^2/ef;b,c,d,a/e,a/f,q^{1-n},q^{-n};q^2,p).
\end{multline*}
\end{theorem}
For $p=0$ this is \cite[Eq. (5.15)]{GR90b} (with corrected misprint).
\begin{theorem}\label{Etrafo4}
For $bc=a^2q^{n+1}$ and $de=aq^{n+1}$,
\begin{multline*}
\sum_{k=0}^{\lfloor n/2 \rfloor} \frac{E(aq^{4k})}{E(a)}
\frac{(b,c;q^3,p)_k}{(aq/b,aq/c;q,p)_k}
\frac{(q^{-n};q,p)_{2k}}{(aq^{n+1};q,p)_{2k}}
\frac{(d,e;q,p)_k}{(aq^3/d,aq^3/e;q^3,p)_k}\: q^k \\
=\frac{(aq;q,p)_n(aq^{2-n}/b;q^3,p)_n}
{(aq/b;q,p)_n(aq^{2-n};q^3,p)_n}
{_{10}\omega_9}(a^2/de;b,c,a/d,a/e,q^{2-n},q^{1-n},q^{-n};q^3,p).
\end{multline*}
\end{theorem}
When $p=0$ this can be recognized as
\cite[Eq. (3.19), $c\to bq^{-n-1}$]{GR90b}.
\begin{theorem}\label{Etrafo5}
For $bcd=a^2q$ and $de=aq^{n+1}$,
\begin{multline*}
\sum_{k=0}^n \frac{E(aq^{4k})}{E(a)}
\frac{(b,c;q^3,p)_k}{(aq/b,aq/c;q,p)_k}
\frac{(d;q,p)_{2k}}{(aq/d;q,p)_{2k}}
\frac{(e,q^{-n};q,p)_k}{(aq^3/e,aq^{n+3};q^3,p)_k}\: q^k \\
=\begin{cases} \displaystyle
f_n \;
{_{10}\omega_9}(a^2/deq;a/dq,a/e,b,c,d,q^{1-n},q^{-n};q^3,p)
& n\not\equiv 2\pmod{3} \\ \displaystyle
g_n \;
{_{10}\omega_9}(a^2/de;a/d,a/e,b,c,dq,q^{2-n},q^{-n};q^3,p)
& n\not\equiv 1\pmod{3} \\ \displaystyle
h_n \;
{_{10}\omega_9}(a^2q/de;aq/d,a/e,b,c,dq^2,q^{2-n},q^{1-n};q^3,p)
& n\not\equiv 0\pmod{3},
\end{cases}
\end{multline*}
with
\begin{align*}
f_n&=
\frac{(aq^{3-\sigma},aq^{3-\sigma}/bc,
aq^{3-\sigma}/bd,aq^{3-\sigma}/cd;q^3,p)_{(n+\sigma)/3}}
{(aq^{3-\sigma}/b,aq^{3-\sigma}/c,
aq^{3-\sigma}/d,aq^{3-\sigma}/bcd;q^3,p)_{(n+\sigma)/3}}  \\
g_n&=
\frac{(aq^{3-\sigma},aq^{3-\sigma}/bc,
aq^{2-\sigma}/bd,aq^{2-\sigma}/cd;q^3,p)_{(n+\sigma)/3}}
{(aq^{3-\sigma}/b,aq^{3-\sigma}/c,
aq^{2-\sigma}/d,aq^{2-\sigma}/bcd;q^3,p)_{(n+\sigma)/3}}  \\
h_n&=\frac{E(aq^\sigma)E(aq^\sigma/bc)E(aq/bd)E(aq/cd)}
{E(aq^\sigma/b)E(aq^\sigma/c)E(aq/d)E(aq/bcd)} \\
& \qquad\times
\frac{(aq^{3-\sigma},aq^{3-\sigma}/bc,
aq^{1-\sigma}/bd,aq^{1-\sigma}/cd;q^3,p)_{(n+\sigma)/3}}
{(aq^{3-\sigma}/b,aq^{3-\sigma}/c,
aq^{1-\sigma}/d,aq^{1-\sigma}/bcd;q^3,p)_{(n+\sigma)/3}}
\end{align*}
where $\sigma\in\{0,1,2\}$ is fixed by
$n+\sigma\equiv 0\pmod{3}$.
\end{theorem}

\begin{proof}
Using identity \eqref{periodicity3} it readily follows that
both the left- and right-hand sides of the identities
in Theorems~\ref{Etrafo3}--\ref{Etrafo5}, viewed as functions of the
variable $b$, satisfy the periodicity $f(pb)=f(b)$.
If we define $h(b)=\text{LHS}(b)/\text{RHS}(b)-1$
then $h$ is a meromorphic function in $0<|b|<\infty$ with 
that same periodicity and with a finite number of poles in
a period annulus. Such a function is either a constant or has an
equal number of zeros and poles in a period annulus
(poles or zeros of order $j$ counted $j$ times).
So if we can show that within a
period annulus $h(b)=0$ for an infinite number of $b$ then
$h$ must be identically zero.
Without loss of generality we may assume that $q^{m_1}\neq p^{m_2}$
for $m_1$ and $m_2$ positive integers.
It is then enough to show that the 
identities in Theorems~\ref{Etrafo3}--\ref{Etrafo5} hold for $b=q^{-m}$
where $m$ runs over an infinite subset of the integers.
First consider Theorem \ref{Etrafo3}.
For $b=q^{-2m}$ ($m$ a nonnegative integer) it holds as can be seen
by making the  simultaneous replacements
$n\to m, b\to e, c\to f$, $d\to q^{-n}, e\to c$, $f\to d$
in the $g=a/b$ case of Theorem~\ref{Etrafo} and using \eqref{el}.
Theorem \ref{Etrafo4} for $b=q^{-3m}$ ($m$ a nonnegative integer)
follows by making the simultaneous replacements
$n\to m, b\to d, c\to e$, $d\to q^{-n}, e\to c$
in the $f=a/b$ case of Theorem~\ref{Etrafo2} and using \eqref{el}.
Finally we show that the first ($n\not\equiv 2 \pmod{3}$) of the
identities of Theorem~\ref{Etrafo5} holds for $b=q^{-3m}$
($m$ a nonnegative integer). The other two follow in similar manner.
Take Theorem~\ref{Etrafo2} with $f=a/b$ and make the simultaneous replacements
$n\to m$, $b\to e$, $c\to q^{-n}$, $e\to c$.
To the thus obtained identity apply the elliptic Bailey transformation 
\eqref{E109} with $a\to a^2 q^n/e$, $b\to aq^n$, $c\to aq^{n+3}/e$, $d\to dq$,
$e\to c$, $f\to d$ and $g\to a/e$ (so that $\lambda=a^2/deq$).
The result is the first identity of the theorem with $b=q^{-3m}$.
\end{proof}

By appropriately specializing the transformations in the Theorems
\ref{Etrafo3}--\ref{Etrafo5} we obtain several further summations. 
Taking $f=a$ in Theorem~\ref{Etrafo3} leads to an elliptic extension of
\cite[Eq. (5.15), $b\to q^{-n}$]{Gasper89} and
\cite[Eq. (1.9), $c\to q^{-n}$]{Rahman93}.
\begin{corollary}
For $bcd=a^2q$ and $e=q^{n+1}$,
\begin{align*}
\sum_{k=0}^n\frac{E(aq^{3k})}{E(a)} &
\frac{(b,c,d;q^2,p)_k}{(aq/b,aq/c,aq/d;q,p)_k}
\frac{(a,e,q^{-n};q,p)_k}{(q^2,aq^2/e,aq^{n+2};q^2,p)_k}\: q^k \\
&=\frac{(aq,aq/bc;q,p)_n(aq^{1-n}/b,aq^{1-n}/c;q^2,p)_n}
{(aq/b,aq/c;q,p)_n(aq^{1-n},aq^{1-n}/bc;q^2,p)_n}\\
&=\frac{(aq^{2-\sigma},aq^{2-\sigma}/bc,aq^{2-\sigma}/bd,
aq^{2-\sigma}/cd;q^2,p)_{(n+\sigma)/2}}
{(aq^{2-\sigma}/b,aq^{2-\sigma}/c,aq^{2-\sigma}/d,
aq^{2-\sigma}/bcd;q^2,p)_{(n+\sigma)/2}},
\end{align*}
where $\sigma\in\{0,1\}$ is determined by $n+\sigma\equiv 0\pmod{2}$.
\end{corollary}
To obtain a summation formula by choosing $d=a$ in Theorem \ref{Etrafo3}
requires a bit more work. First observe that by this choice for $d$ 
we have $bc=aq$ so that
\begin{equation*}
\frac{(aq/bc;q,p)_n}{(aq^{1-n}/bc;q^2,p)_n}=
\begin{cases}
\displaystyle \frac{(q;q^2,p)_{n/2}}{(q^{-n};q^2,p)_{n/2}} & 
\text{for $n$ even} \\[2mm]
0 & \text{for $n$ odd.}
\end{cases}
\end{equation*}
Next note that
${_{10}\omega_9}(a^2/ef;a,b,c,a/e,a/f,q^{1-n},q^{-n};q^2,p)$
with $ef=aq^{n+1}$ and $bc=aq$ reduces to 
${_8\omega_7}(a^2/ef;b,c,a/e,a/f,q^{-n};q^2,p)$,
which for $n$ even can be summed by \eqref{e87} to give
\begin{equation*}
\frac{(q^2,q/a,bq/e,aq^2/be;q^2,p)_{n/2}}{(q/e,aq^2/e,q^2/b,bq/a;q^2,p)_{n/2}}.
\end{equation*}
Combining all of the above and using \eqref{el} we get the following 
generalization of summation \cite[Eq. (6.14)]{GS83} by Gessel and Stanton.
\begin{corollary}\label{EGS}
For $bc=aq$ and $de=aq^{n+1}$,
\begin{multline*}
\sum_{k=0}^n \frac{E(aq^{3k})}{E(a)}
\frac{(a,b,c;q^2,p)_k}{(q,aq/b,aq/c;q,p)_k}
\frac{(d,e,q^{-n};q,p)_k}{(aq^2/d,aq^2/e,aq^{n+2};q^2,p)_k}\: q^k \\
=\begin{cases}\displaystyle
\frac{(aq^2,aq^2/bc,aq^2/bd,aq^2/cd;q^2,p)_{n/2}}
{(aq^2/b,aq^2/c,aq^2/d,aq^2/bcd;q^2,p)_{n/2}} & \text{for $n$ even} \\[3mm]
0 & \text{for $n$ odd.} 
\end{cases}
\end{multline*}
\end{corollary}

In exactly the same manner we can derive two summations from Theorem
\ref{Etrafo4}. First when $e=a$ we get the elliptic analogue of
\cite[Eq. (5.22), $b\to q^{n+1}$]{Gasper89}.
\begin{corollary}
For $bc=a^2q^{n+1}$ and $d=q^{n+1}$,
\begin{multline*}
\sum_{k=0}^{\lfloor n/2 \rfloor} \frac{E(aq^{4k})}{E(a)}
\frac{(b,c;q^3,p)_k}{(aq/b,aq/c;q,p)_k}
\frac{(q^{-n};q,p)_{2k}}{(aq^{n+1};q,p)_{2k}}
\frac{(a,d;q,p)_k}{(aq^3/d,q^3;q^3,p)_k}\: q^k \\
=\frac{(aq;q,p)_n(aq^{2-n}/b;q^3,p)_n}
{(aq/b;q,p)_n(aq^{2-n};q^3,p)_n}.
\end{multline*}
\end{corollary}

The second summation follows from $c=a$. Then $b=aq^{n+1}$ and 
\begin{equation*}
\frac{(aq^{2-n}/b;q^3,p)_n}{(aq/b;q,p)_n}=
\frac{(q^{1-2n};q^3,p)_n}{(q^{-n};q,p)_n}=
0 \qquad \text{for $n\equiv 2 \pmod{3}$.}
\end{equation*}
Further observe that
${_{10}\omega_9}(a^2/de;a,b,a/d,a/e,q^{2-n},q^{1-n},q^{-n};q^3,p)$
with $b=aq^{n+1}$ and $de=aq^{n+1}$ reduces to
${_8\omega_7}(a^2/de;a/d,a/e,b,q^{-n},q^{1-n};q^3,p)$,
which for $n\not\equiv 2\pmod{3}$ can be summed to give
\begin{equation*}
\frac{(aq^{2-n},q^3,dq^{1-2n}/a,q^{2-n}/d;q^3,p)_{\lfloor n/3 \rfloor}}
{(q^{2-n}/a,q^{1-2n},aq^3/d,dq^{2-n};q^3,p)_{\lfloor n/3 \rfloor}}.
\end{equation*}
After a few simplification we arrive at the following summation.
\begin{corollary}\label{C2}
For $b=aq^{n+1}$ and $cd=aq^{n+1}$,
\begin{multline*}
\sum_{k=0}^{\lfloor n/2 \rfloor} \frac{E(aq^{4k})}{E(a)}
\frac{(a,b;q^3,p)_k}{(q,aq/b;q,p)_k}
\frac{(q^{-n};q,p)_{2k}}{(aq^{n+1};q,p)_{2k}}
\frac{(c,d;q,p)_k}{(aq^3/c,aq^3/d;q^3,p)_k}\: q^k \\
=\begin{cases}\displaystyle
\frac{(aq^3,aq^3/bc,aq^3/bd;q^3,p)_{\lfloor n/3 \rfloor}}
{(aq^3/c,aq^3/d,aq^3/bcd;q^3,p)_{\lfloor n/3 \rfloor}}
& n\not\equiv 2\pmod{3} \\[3mm]
0 & n\equiv 2\pmod{3}.
\end{cases}
\end{multline*}
\end{corollary}
Finally we turn to Theorem~\ref{Etrafo5}.
The choice $e=a$ immediately gives an elliptic analogue of
\cite[Eq. (5.22), $b\to q^{-n}$]{Gasper89}.
\begin{corollary}
For $bcd=a^2q$ and $d=q^{n+1}$,
\begin{multline*}
\sum_{k=0}^n \frac{E(aq^{4k})}{E(a)}
\frac{(b,c;q^3,p)_k}{(aq/b,aq/c;q,p)_k}
\frac{(d;q,p)_{2k}}{(aq/d;q,p)_{2k}}
\frac{(a,q^{-n};q,p)_k}{(q^3,aq^{n+3};q^3,p)_k}\: q^k \\
=\begin{cases} \displaystyle
\frac{(aq^3,aq^3/bc,aq^3/bd,aq^3/cd;q^3,p)_{n/3}} 
{(aq^3/b,aq^3/c,aq^3/d,aq^3/bcd;q^3,p)_{n/3}} 
& n\equiv 0\pmod{3} \\[3mm] \displaystyle
\frac{(aq,aq/bc,aq/bd,aq/cd;q^3,p)_{(n+2)/3}} 
{(aq/b,aq/c,aq/d,aq/bcd;q^3,p)_{(n+2)/3}} 
& n\equiv 1\pmod{3} \\[3mm] \displaystyle
\frac{(aq^2,aq^2/bc,aq/bd,aq/cd;q^3,p)_{(n+1)/3}} 
{(aq^2/b,aq^2/c,aq/d,aq/bcd;q^3,p)_{(n+1)/3}} 
& n\equiv 2\pmod{3}.
\end{cases}
\end{multline*}
\end{corollary}
If, on the other hand, we set $c=a$ in Theorem~\ref{Etrafo5}
and perform a calculation similar to the one employed in the derivation
of Corollaries~\ref{EGS} and \ref{C2} we get the elliptic extension of
\cite[Eq. (4.6d)]{Chu95}.
\begin{corollary}\label{CorChu}
For $bc=aq$ and $cd=aq^{n+1}$,
\begin{multline*}
\sum_{k=0}^n \frac{E(aq^{4k})}{E(a)}
\frac{(a,b;q^3,p)_k}{(q,aq/b;q,p)_k}
\frac{(c;q,p)_{2k}}{(aq/c;q,p)_{2k}}
\frac{(d,q^{-n};q,p)_k}{(aq^3/d,aq^{n+3};q^3,p)_k}\: q^k \\
=\begin{cases}\displaystyle
\frac{(q,q^2,aq^3,b^2/a;q^3,p)_{n/3}}{(bq,bq^2,b/a,aq^3/b;q^3,p)_{n/3}}
& n\equiv 0\pmod{3} \\[3mm]
0 &  n\not\equiv 0 \pmod{3}.
\end{cases}
\end{multline*}
\end{corollary}
Similarly, taking $d=a$ in Theorem~\ref{Etrafo5} yields the last summation
of this section.
\begin{corollary}
For $bc=aq$ and $d=q^{n+1}$,
\begin{multline*}
\sum_{k=0}^n \frac{E(aq^{4k})}{E(a)}
\frac{(b,c;q^3,p)_k}{(aq/b,aq/c;q,p)_k}
\frac{(a;q,p)_{2k}}{(q;q,p)_{2k}}
\frac{(d,q^{-n};q,p)_k}{(aq^3/d,aq^{n+3};q^3,p)_k}\: q^k \\
=\begin{cases} \displaystyle
\frac{(aq^3,q^2/b,q^2/c;q^3,p)_{n/3}}{(q^2/bc,aq^3/b,aq^3/c;q^3,p)_{n/3}} 
& n\equiv 0\pmod{3} \\[3mm]
0 & n\equiv 1\pmod{3} \\
\displaystyle
\frac{(aq^2,q/b,q/c;q^3,p)_{(n+2)/3}}{(q/bc,aq^2/b,aq^2/c;q^3,p)_{(n+2)/3}} 
& n\equiv 2\pmod{3}.
\end{cases}
\end{multline*}
\end{corollary}
For $p=0$ this is \cite[Eq. (3.21), $a\to q^{-n}$]{GR90b}.

Before concluding this section let us remark that there
is a nice corollary to Corollary~\ref{cor1}
in the form of a nontrivial determinant evaluation.
For $p=0$ this was observed by Andrews and Stanton
and the following theorem is a direct elliptic extension to their
\cite[Thm. 8]{AS98}. (In fact Andrews and Stanton did not use the
$p=0$ version of Corollary~\ref{cor1} but of the $p=0$, $n$ odd
instance of Corollary~\ref{EGS} in their proof.)
\begin{theorem}
For $x,y$ indeterminates and $n$ a positive integer there holds
\begin{multline*}
\det_{1\leq i,j\leq n}\biggl(
\frac{(y q^{1-i}/x,q^{2-i}/xy,q^{2-4i}/x^2;q^2,p)_{i-j}}
{(q^{2-2i}/xy,yq^{1-2i}/x,q^{i+1};q,p)_{i-j}}\biggr)\\=
\prod_{i=1}^n
\frac{(q,x^2q^{2i-2};q,p)_i}{(q,x^2q^{2i-2};q^2,p)_i}
\frac{(xyq^{i-1},xq^i/y;q^2,p)_i}{(xyq^{i-1},xq^i/y;q,p)_i}.
\end{multline*}
\end{theorem}
Let us note that for $i<j$ we need the elliptic analogue of the $q$-shifted
factorial \eqref{EP} with negative subscript. Hence
$1/(q^{i+1};q,p)_{i-j}=(q^{2i-j+1};q,p)_{j-i}=0$ for $2i-j+1\leq 0$.

\begin{proof}
To prove the theorem we establish the Gauss decomposition or ``LU''
factorization \cite{Krattenthaler99} of the matrix $M$ featuring in the 
determinant.
Let $M_n=(M_{i,j})_{1\leq i,j\leq n}$ with
\begin{equation*}
M_{i,j}=\frac{(y q^{1-i}/x,q^{2-i}/xy,q^{2-4i}/x^2;q^2,p)_{i-j}}
{(q^{2-2i}/xy,yq^{1-2i}/x,q^{i+1};q,p)_{i-j}},
\end{equation*}
and
$U_n=(U_{i,j})_{1\leq i,j\leq n}$, with
\begin{equation*}
U_{i,j}=
(-1)^{i+j} q^{(i-j)(i+j-7)/2}
\frac{E(x^2q^{3i-2})(q^i;q,p)_{2j-2i}(q^{3-3j}/x^2;q,p)_{j-i}}
{E(x^2q^{i+2j-2})(q^2,q^{4-4j}/x^2,q^{3-2j}/x^2;q^2,p)_{j-i}} 
\end{equation*}
for $i\leq j$, and
\begin{equation*}
U_{i,j}=0
\end{equation*}
for $i>j$.
Next we calculate the product of the above two matrices
\begin{multline*}
(M_n\cdot U_n)_{i,j}=
\sum_{k=1}^j M_{i,k} U_{k,j} \\
=\frac{(y q^{1-i}/x,q^{2-i}/xy,q^{2-4i}/x^2;q^2,p)_{i-1}
(q^j,q^{1-j},q^{3-3j}/x^2;q,p)_{j-1}q^{j-1}}
{(q^{2-2i}/xy,yq^{1-2i}/x,q^{i+1};q,p)_{i-1}
(q^{4-4j}/x^2,q^{1-2j}/x^2,q^2;q^2,p)_{j-1}}  \\
\times \sum_{k=0}^{j-1}
\frac{E(x^2q^{3k+1})}{E(x^2q)}
\frac{(xq^{i+1}/y,xyq^i,q^{1-2i};q,p)_k(x^2q,x^2q^{2j},q^{2-2j};q^2,p)_k q^k}
{(xyq^{2-i},xq^{3-i}/y,x^2q^{2i+2};q^2,p)_k(q,q^{2-2j},x^2q^{2j};q,p)_k}.
\end{multline*}
To proceed we observe that the sum over $k$ can be carried out by
the $e=a$ case of Corollary~\ref{cor1} so that the last line in the
above equation may be replaced by
\begin{equation*}
\frac{(x^2 q^3,xq^{i+2}/y,xyq^{i+1},q^{2-2i};q^2,p)_{j-1}}
{(q,xyq^{2-i},xq^{3-i}/y,x^2q^{2i+2};q^2,p)_{j-1}}.
\end{equation*}
We learn two things from this result. First, that the matrix $L_n=
(L_{i,j})_{1\leq i,j\leq n}=M_n\cdot U_n$
is lower triangular (since $(q^{2-2i};q^2,p)_{j-1}=0$ for $j>i$)
and, second, that its diagonal entries are given by
\begin{equation*}
L_{i,i}=\frac{(q^2,x^2q^{2i-1};q,p)_{i-1}
(xyq^{i+1},xq^{i+2}/y;q^2,p)_{i-1}}
{(q^3,x^2q^{2i};q^2,p)_{i-1}(xyq^i,xq^{i+1}/y;q,p)_{i-1}}.
\end{equation*} 
The calculation of $\det(M_n)$ is now done; by $M_n\cdot U_n=L_n$
we get $\det(M_n)\det(U_n)=\det(L_n)$, but $\det(U_n)=1$ by the fact
that $U_n$ is an upper-triangular matrix with $1$'s along the diagonal.
Hence we only need to compute the determinant of $L_n$ which is the
product of its diagonal entries, resulting in the right-hand side of
the theorem.
\end{proof}

\section{An elliptic C$_n$ Jackson sum}\label{secC}
Building on earlier work in \cite{GK97}, Schlosser
proved a multidimensional extension of Jackson's
$_8\phi_7$ summation~\cite{Schlosser00}.
Here we show that by a generalization of a determinant lemma of 
Krattenthaler (see Lemmas \ref{lemK} and \ref{lemeldet} below)
Schlosser's C$_n$ Jackson sum can readily be generalized to the
elliptic case. This is the content of our next theorem.
\begin{theorem}\label{thmEcn}
For $x_1,\dots,x_n$, $a,b,c,d$ and $e$ indeterminates and $N$ a nonnegative
integer such that $a^2q^{N-n+2}=bcde$ there holds
\begin{multline}\label{Jackson87sum}
\sum_{k_1,\dots,k_n=0}^N
\prod_{1\leq i<j\leq n}\biggl(\frac{E(q^{k_i-k_j}x_i/x_j)}{E(x_i/x_j)}
\frac{E(ax_ix_jq^{k_i+k_j})}{E(ax_ix_jq^N)}\biggr) \\
\times\prod_{i=1}^n \frac{E(ax_i^2q^{2k_i})}{E(ax_i^2)}
\frac{(ax_i^2,bx_i,cx_i,dx_i,ex_i,q^{-N};q,p)_{k_i}q^{ik_i}}
{(q,aqx_i/b,aqx_i/c,aqx_i/d,aqx_i/e,ax_i^2q^{N+1};q,p)_{k_i}}\\
=\prod_{i=1}^n\frac{(aqx_i^2,aq^{2-i}/bc,aq^{2-i}/bd,aq^{2-i}/cd;q,p)_N}
{(aq^{2-n}/bcdx_i,aqx_i/b,aqx_i/c,aqx_i/d;q,p)_N}.
\end{multline}
\end{theorem}
As remarked above, to prove this result we need the elliptic analogue of the
following determinant lemma due to 
Krattenthaler~\cite[Lemma 34]{Krattenthaler95} 
(see also \cite[Lemma 5]{Krattenthaler99}), which was crucial in the
proof of the $p=0$ case of \eqref{Jackson87sum}~\cite{Schlosser00}.
\begin{lemma}\label{lemK}
Let $X_1,\dots,X_n,A_2,\dots,A_n$ and $C$ be indeterminates.
If, for $j=0,\dots,n-1$, $P_j$ is a Laurent polynomial of degree less than 
or equal to $j$ such that $P_j(C/X)=P_j(X)$, then
\begin{multline*}
\det_{1\leq i,j\leq n}\Bigl(P_{j-1}(X_i)
\prod_{k=j+1}^n(1-A_kX_i)(1-CA_k/X_i)\Bigr) \\
=\prod_{1\leq i<j\leq n}A_j X_j(1-X_i/X_j)(1-C/X_iX_j)
\prod_{i=1}^n P_{i-1}(1/A_i).
\end{multline*}
\end{lemma}
Here the degree of a Laurent polynomial $P(x)=\sum_{i=M}^N a_i x^i$ with
$a_N\neq 0$ is defined to be $N$,
and the empty product $\prod_{k=j+1}^n (1-A_kX_i)(1-CA_k/X_i)$ for
$j=n$ is defined to be $1$.
For a proof of this lemma we refer to \cite{Krattenthaler95}.

The needed elliptic analogue of the previous lemma can be stated as follows.
\begin{lemma}\label{lemeldet}
Let $X_1,\dots,X_n,A_2,\dots,A_n$ and $C$ be indeterminates
and $E$ the elliptic function defined in \eqref{defE}.
If, for $j=0,\dots,n-1$, $P_j$ is analytic in $0<|x|<\infty$ with
periodicity $P_j(px)=(C/x^2p)^j P_j(x)$ and symmetry $P_j(C/x)=P_j(x)$, then
\begin{multline}\label{eldet}
\det_{1\leq i,j\leq n}\Bigl(P_{j-1}(X_i)
\prod_{k=j+1}^n E(A_kX_i)E(CA_k/X_i)\Bigr) \\
=\prod_{1\leq i<j\leq n}A_j X_jE(X_i/X_j)E(C/X_iX_j)
\prod_{i=1}^n P_{i-1}(1/A_i).
\end{multline}
\end{lemma}

\begin{proof}
View both sides of \eqref{eldet} as a function of the variable $X_i$ 
$(i=1,\dots,n)$, and write $L(X_i)$ ($R(X_i)$) for the left(right)-hand 
side. From the periodicity property \eqref{periodicity1}
and the periodicity of $P_j$, we find that 
$$F(X_i)=(pX_i^2/C)^{n-1}F(pX_i),$$
where $F=L,R$. 
As a result the function $f$, defined as the ratio
of $L$ over $R$, satisfies the periodicity $f(X_i)=f(pX_i)$.
Since $E(x)$ and $P_j(x)$ are analytic in $0<|x|<\infty$,
the only possible poles of $f$ are the zero's of $R$.
Since $E(x)$ has simple zeros at
$x=p^k$ ($k\in\Z$), the zeros of $R$ are
$X_i=p^k X_j$ and $X_i=p^k C/X_j$
where $k\in\Z$ and $j=1,\dots,i-1,i+1,\dots,n$.
First consider $X_i=p^k X_j$. When inserted into
the determinant it follows from \eqref{periodicity1}
and 
\begin{equation}\label{periodicity2}
P_j(x)=(x^{2} p^{k}/C)^{jk} P_j(xp^k),
\end{equation}
that the $i$-th and $j$-th row become proportional
(with proportionality constant $(Cp^{-k}/X_j^2)^{k(n-1)}$).
Next, when $X_i=p^k C/X_j$ it follows from
\eqref{periodicity1}, \eqref{periodicity2} and the symmetry
$P_j(C/x)=P_j(x)$ that the 
$i$-th and $j$-th row once again become proportional
(with proportionality constant $(X_j^2p^{-k}/C)^{k(n-1)}$).
We may therefore conclude that 
$L$ vanishes at the zeros of $R$, so that, according to 
Liouville's theorem, $f$ must be constant.

To conclude the proof we only need to show the validity of 
\eqref{eldet} for some appropriately chosen values of $X_1,\dots,X_n$.
A good choice is 
\begin{equation}\label{Xspec}
X_i=1/A_i, \quad \text{$i=1,\dots,n$}.
\end{equation}
Since $\prod_{k=j+1}^n E(X_i/A_k)=0$ for $j<i$, this
leaves the determinant of an upper-triangular matrix which 
evaluates to
\begin{equation*}
\prod_{i=1}^n \prod_{j=i+1}^n E(A_j/A_i)E(CA_iA_j)P_{i-1}(1/A_i).
\end{equation*}
Clearly this corresponds to the right-hand side of \eqref{eldet} under
the specialization \eqref{Xspec}, and we are done.
\end{proof}

Choosing $A_i=Aq^{n-i}$ and $P_i(X)=(BXq^{n-i-1},BCq^{n-i-1}/X;q,p)_i$ 
in \eqref{eldet}, and using \eqref{el} and $\sum_{j=1}^n(j-1)(n-j)
=\binom{n}{3}$, we obtain the following nice
corollary of Lemma~\ref{lemeldet}.
\begin{corollary}\label{determ}
For $X_1,\dots,X_n,A,B$ and $C$ indeterminates,
\begin{multline*}
\det_{1\leq i,j\leq n}\biggl(
\frac{(AX_i,AC/X_i;q,p)_{n-j}}
{(BX_i,BC/X_i;q,p)_{n-j}}\biggr) \\ =
A^{\binom{n}{2}}q^{\binom{n}{3}} 
\prod_{1\leq i<j\leq n}X_j E(X_i/X_j)E(C/X_iX_j)
\prod_{i=1}^n
\frac{(B/A,ABCq^{2n-2i};q,p)_{i-1}}
{(BX_i,BC/X_i;q,p)_{n-1}}.
\end{multline*}
\end{corollary}
We remark that this can be written as the following determinant identity
for theta functions:
\begin{multline*}
\det_{1\leq i,j\leq n}\Bigl(T_{n-j}(A+X_i)T_{n-j}(A+C-X_i) \\
\times T_{j-1}(B+X_i+n-j)T_{j-1}(B+C+n-j-X_i)\Bigr) \\
=\prod_{1\leq i<j\leq n}\vartheta_1(X_i-X_j)\vartheta_1(C-X_i-X_j)
\prod_{i=1}^n T_{i-1}(B-A)T_{i-1}(A+B+C+2n-2i),
\end{multline*}
where $T_n(x)=\prod_{k=0}^{n-1}\vartheta_1(x+k)$ and $\vartheta_1(x)$
a standard theta function \cite{WW96}, 
\begin{equation*}
\vartheta_1(x)=2\sum_{k=0}^{\infty}(-1)^n p^{(2n+1)^2/4}\sin(2n+1)x
=ip^{1/4}e^{-ix}(p^2;p^2)_{\infty}E(e^{2ix};p^2).
\end{equation*}
For $n=2$ this is nothing but the well-known identity
\begin{multline*}
\vartheta_1(u+x)\vartheta_1(u-x)\vartheta_1(v+y)\vartheta_1(v-y)-
\vartheta_1(u+y)\vartheta_1(u-y)\vartheta_1(v+x)\vartheta_1(v-x)\\=
\vartheta_1(x+y)\vartheta_1(x-y)\vartheta_1(u+v)\vartheta_1(u-v),
\end{multline*}
equivalent to \eqref{Esum}.

\begin{proof}[Proof of Theorem~\ref{thmEcn}]
By Corollary~\ref{determ} with $X_i\to q^{-k_i}/x_i$ and $C\to a$, and
$E(x)=-xE(1/x)$ we can trade the double product 
$$\prod_{1\leq i<j\leq n}E(q^{k_i-k_j}x_i/x_j)E(ax_ix_jq^{k_i+k_j})$$
for a determinant. If we also choose 
$B=q^{2-n}/c$ and $A=b/a$, and use \eqref{el}
and $\sum_{j=1}^n\binom{j-1}{2}=\binom{n}{3}$,
the left-hand side of \eqref{Jackson87sum}
can be rewritten as
\begin{multline}\label{s1}
q^{-3\binom{n}{3}} 
\prod_{1\leq i<j\leq n}\biggl(
\frac{a^2x_j/bc^2}{E(x_j/x_i)E(ax_ix_jq^N)}\biggr)
\prod_{i=1}^n\frac{1}{(aq^{2-n}/bc,bq^{n-2i+2}/c;q,p)_{i-1}} \\
\times\det_{1\leq i,j\leq n}\biggl(
(cx_i,c/ax_i;q,p)_{j-1}
(bx_i,b/ax_i;q,p)_{n-j} \\
\times{_8\omega_7}(ax_i^2;bx_iq^{n-j},cx_iq^{j-1},dx_i,ex_i,q^{-N};q,p)
\biggr).
\end{multline}
Applying the elliptic $_8\omega_7$ sum of Theorem~\ref{e87} and again 
using \eqref{el} as well as
$\sum_{j=1}^n(n-j)(n+j-3)=4\binom{n}{3}$, we arrive at the following
expression for the left-hand side of \eqref{Jackson87sum}
\begin{multline*}
q^{-\binom{n}{3}}
\prod_{1\leq i<j\leq n}\biggl(
\frac{ax_jq^N/b}{E(x_j/x_i)E(ax_ix_jq^N)}\biggr)
\prod_{i=1}^n \frac{(q^{2-n}/cx_i,ax_iq^{N-n+2}/c;q,p)_{n-1}}
{(bq^{n-2i+2}/c,aq^{N-n+2}/bc;q,p)_{i-1}} \\
\times \prod_{i=1}^n
\frac{(ax_i^2q,aq^{2-i}/bc,aq^{2-i}/bd,aq^{2-i}/cd;q,p)_N}
{(aqx_i/b,aqx_i/c,ax_iq/d,aq^{2-n}/bcdx_i;q,p)_N} \\
\times\det_{1\leq i,j\leq n}\biggl(
\frac{(bx_i,bq^{-N}/ax_i;q,p)_{n-j}}
{(q^{2-n}/cx_i,ax_iq^{N-n+2}/c;q,p)_{n-j}}\biggr).
\end{multline*}
By Lemma~\ref{determ} with $X_i\to 1/x_i$, $A\to bq^{-N}/a$, $B\to q^{2-n}/c$
and $C\to aq^N$ the first and third line are found to be reciprocal,
thus resulting in the right-hand side of \eqref{Jackson87sum}.
\end{proof}

\section{Discussion}
Of course the summations and transformations obtained in this paper
for elliptic hypergeometric series
are only a tip of the iceberg. Many more results
for terminating, balanced, very-well-poised, basic hypergeometric series
admit elliptic generalizations.
In particular all the multivariable balanced, very-well-poised
summation and transformation theorems of
\cite{BM97,BS98,DG92,Milne88,Milne89,Milne94,ML95,MN96,Schlosser97}
should admit elliptic counterparts.
However, the methods of proof applied in these papers does not
simply carry over the the elliptic case. In particular,
the multivariable Jackson sums (from which most of the other results
can be derived in ways well-tailored for elliptic generalization)
are usually proved using
simpler identities for series that are not both balanced and very-well poised.
This is unlike the one-dimensional Jackson sum which can be
proved simply by induction, without relying on other results --
a method of proof that readily carries over to the elliptic case.
Indeed the only higher-dimensional elliptic Jackson sum that we were able to
prove so far is the one stated in Theorem~\ref{thmEcn}.

One might expect that at least ``the corresponding 
$_{10}\omega_9$ transformation'' should be accessible with the
techniques presented in this paper. However, all our attempts to 
find a C$_n$ $_{10}\omega_9$ transformation that implies
Theorem~\ref{thmEcn} failed dismally. Surprisingly though,
our failed attempts did suggest how to somewhat change
Theorem~\ref{thmEcn} so that it does admit a generalization
to a transformation. Since to the best of our knowledge
this transformation (in the $p\to 0$ limit) does not appear in the above list
of references we state it here as a conjecture.

First we we need some more notation.
Following Macdonald's book~\cite{Macdonald95} we set
$$|\lambda|=\sum_{i\geq 1}\lambda_i \qquad \text{and} \qquad
n(\lambda)=\sum_{i\geq 1} (i-1)\lambda$$
for $\lambda$ a partition (i.e., $\lambda=(\lambda_1,\lambda_2,\dots)$
with $\lambda_i\geq \lambda_{i+1}$ and finitely many $\lambda_i$ nonzero).
For a partition $\lambda$ of exactly $n$ parts (some of which may be zero)
define
\begin{equation*}
(a;q,p)_{\lambda}=\prod_{i=1}^n (ax^{1-j};q,p)_{\lambda_i}
\end{equation*}
and employ the usual condensed notation
\begin{equation*}
(a_1,\dots,a_m;q,p)_{\lambda}=(a_1;q,p)_{\lambda}\dots (a_m;q,p)_{\lambda}.
\end{equation*}
With these preliminaries we define
a C$_n$ analogue of the balanced, very-well-poised, elliptic
hypergeometric series \eqref{omega} by
\begin{multline*}
{_{r+1}\Omega_r}(a_1;a_4,\dots,a_{r+1};q,p) \\
=\sum_{\lambda_1\geq\lambda_2\geq\dots\geq\lambda_n\geq 0} \;
\prod_{i=1}^n\Biggl(
\frac{E(a_1x^{2(1-i)}q^{2\lambda_i})}{E(a_1x^{2(1-i)})}\Biggl)
\frac{(a_1x^{1-n},a_4,\dots,a_{r+1};q,p)_{\lambda}q^{|\lambda|}x^{2n(\lambda)}}
{(qx^{n-1},a_1q/a_4,\dots,a_1q/a_{r+1};q,p)_{\lambda}}\\
\times
\prod_{1\leq i<j\leq n}\biggl(
\frac{E(x^{j-i}q^{\lambda_i-\lambda_j})}{E(x^{j-i})} \;
\frac{E(a_1x^{2-i-j}q^{\lambda_i+\lambda_j})}{E(a_1x^{2-i-j})} \\
\times
\frac{(a_1x^{3-i-j};q,p)_{\lambda_i+\lambda_j}
(x^{j-i+1};q,p)_{\lambda_i-\lambda_j}}
{(a_1qx^{1-i-j};q,p)_{\lambda_i+\lambda_j}
(qx^{j-i-1};q,p)_{\lambda_i-\lambda_j}}\biggr),
\end{multline*}
where $(a_4\dots a_{r+1})^2=a_1^{r-3}q^{r-5}x^{2-2n}$.
For reasons of convergence we again insist that 
one of the $a_i$ $(i=4,\dots,r+1)$
is of the form $q^{-N}$ with $N$ a nonnegative integer,
so that the only nonvanishing contributions
to the above sum come from $\lambda_1\leq N$.
Observe that for $x=1$ the double product in the summand simplifies to 
a multinomial coefficient, i.e., to
$\prod_{1\leq i<j\leq n}(j-i+1-\delta_{\lambda_i,\lambda_j})/(j-i)=
n!/(m_0!m_1!\dots m_N!)$, where $m_k$ is the number of
parts of size $k$ in the partition $\lambda=(\lambda_1,\dots,\lambda_n)$.
Since $\lim_{x\to 1}(a;q,p)_{\lambda}=\prod_{i=1}^n (a;q,p)_{\lambda_i}$ and
\begin{equation*}
\sum_{\lambda}
\frac{n!}{m_0!\dots m_N!}\prod_{i=1}^n a_{\lambda_i}=
\sum_{\substack{0\leq m_0,\dots,m_N\leq n \\ m_1+\cdots+m_N=n}}
\frac{n!}{m_0!\dots m_N!}\prod_{i=0}^N a_i^{m_i}=
\Big(\sum_{i=0}^N a_i\Big)^n,
\end{equation*}
where $\lambda=(\lambda_1,\dots,\lambda_n)=(0^{m_0}1^{m_1}\dots N^{m_N})$),
we may conclude that
$$\lim_{x\to 1}{_{r+1}\Omega_r}(a_1;a_4,\dots,a_{r+1};q,p)
=\bigl({_{r+1}\omega_r}(a_1;a_4,\dots,a_{r+1};q,p)\bigr)^n.$$

Computer assisted experiments suggest the following 
C$_n$ version of the $_{10}\omega_9$ transformation \eqref{E109}.
\begin{conjecture}
Let $bcdefg x^{n-1}=a^3q^{N+2}$ and $\lambda=a^2q/bcd$.
Then
\begin{multline*}
{_{10}\Omega_9}(a;b,c,d,e,f,g,q^{-N};q,p) \\
=\frac{(aq,aq/ef,\lambda q/e,\lambda q/f;q,p)_{(N^n)}}
{(aq/e,aq/f,\lambda q/ef,\lambda q;q,p)_{(N^n)}}\;
{_{10}\Omega_9}(\lambda;\lambda b/a,\lambda c/a,\lambda d/a,e,f,g,q^{-N};q,p).
\end{multline*}
\end{conjecture}
For $cd=aq$ this implies
\begin{corollary}\label{Omega87}
For $bcfg x^{n-1}=a^2q^{N+1}$ there holds
\begin{equation*}
{_8\Omega_7}(a;b,c,d,e,q^{-N};q,p)
=\frac{(aq,aq/bc,aq/bd,aq/cd;q,p)_{(N^n)}}
{(aq/b,aq/c,aq/d,aq/bcd;q,p)_{(N^n)}}.
\end{equation*}
\end{corollary}
As remarked earlier we were unable to trace the $p=0$ case of the above
two results in the literature, but we did find that letting $d$
tend to infinity after setting $p=0$, Corollary \ref{Omega87} reduces
to a multivariable analogue of Rogers' $_6\phi_5$ sum
due to van Diejen \cite[Thm.~3]{vanDiejen97}.

Another challenging problem is to find nontrivial transformations
based on the inverse pair given in \eqref{invp} for all positive integers $r$.
The only result for general $r$ obtained so far in the not-so-deep 
Theorem \ref{thmr}, which we were unable to generalize to a transformation.
The problem with the type of transformations derived in
section~\ref{secST} appears to be that increasing $r$ has the effect of 
decreasing the number of available free parameters. For example, when we mimic
the derivation of Theorems \ref{Etrafo} and \ref{Etrafo2} but choose $r=4$ in
\eqref{invp} we no longer obtain a transformation for a $_{10}\omega_9$,
but the less appealing quartic transformation
\begin{multline*}
\sum_{k=0}^n \frac{E(aq^{5k})}{E(a)}
\frac{(b^2/aq^2;q,p)_k}{(a^2q^6/b^2;q^4,p)_k}
\frac{(aq/b,aq^2/b,aq^3/b;q^2,p)_k}{(b,bq,bq^2;q^3,p)_k}
\frac{(abq^{4n},q^{-4n};q^4,p)_k}{(q^{1-4n}/b,aq^{4n+1};q,p)_k}\: q^k \\
=\frac{(aq;q,p)_{4n}(q^4,b^3/aq^2;q^4,p)_n}
{(b;q,p)_{4n}(ab,a^2q^6/b^2;q^4,p)_n} \\
\times
\sum_{k=0}^n \frac{E(abq^{8k-4})}{E(ab/q^4)}
\frac{(ab/q^4,a^2q^2/b^2,b,b/q,b/q^2,b/q^3;q^4,p)_k}
{(q^4,b^3/aq^2,a,aq,aq^2,aq^3;q^4,p)_k}\: q^{4k},
\end{multline*}
which contains only two indeterminates.
Moreover its counterpart (in the sense of $(q^{-4n};q^4,p)_k \leftrightarrow
(q^{-n};q,p)_k)$) no longer seems to allow for a transformation at all,
admitting just
\begin{multline*}
\sum_{k=0}^n \frac{E(a^2q^{5k})}{E(a^2)}
\frac{(a^2;q^4,p)_k}{(q;q,p)_k}\,
\frac{(a,aq,aq^2;q^3,p)_k}{(a,aq,aq^2;q^2,p)_k}
\frac{(aq^{n+1},q^{-n};q,p)_k}{(aq^{3-n},a^2q^{n+4};q^4,p)_k}\: q^k \\
=\begin{cases}\displaystyle
\frac{(q,q^2,q^3,a^2q^4;q^4,p)_{n/4}}{(aq^2,aq^3,aq^4,q/a;q^4,p)_{n/4}}
& n\equiv 0\pmod{4} \\[3mm]
0 &  n\not\equiv 0 \pmod{0},
\end{cases}
\end{multline*}
which generalizes the quadratic and cubic summations of Corollaries
\ref{EGS} and \ref{CorChu}.

\subsection*{Acknowledgements}
This work is supported by a fellowship of the Royal
Netherlands Academy of Arts and Sciences.

\bibliographystyle{amsplain}

\end{document}